\documentclass[10pt]{amsart}
\usepackage{amsmath}
\usepackage{amscd}
\usepackage{graphics}
\theoremstyle{plain}

\newtheorem{Thm}{Theorem}

\newtheorem{Cor}[Thm]{Corollary}
\newtheorem{Prop}[Thm]{Proposition}

\theoremstyle{definition}
\newtheorem{Def}{Definition}

\theoremstyle{remark}
\newtheorem{Rem}{Remark}

\newtheorem{Rems}{Remarks}

\def\CM{\mathcal M}
\def\CH{\mathcal H}

\def\IR{\mathbb R}
\def\IZ{\mathbb Z}

\def\s2x{\hbox{$S^2 \times S^2$}}
\def\z{\zeta}

\def\a{\alpha}
\def\b{\beta}

\def\e{\eta}
\def\g{\gamma}
\def\S{\Sigma}
\def\bdy{\partial}

\def\vo{\overset{\circ}{\nu}}

\begin{document}
\title[]{Signatures of Lefschetz Fibrations }
\author{Burak Ozbagci}
\address{Department of Mathematics\\University of California\\
Irvine, CA 92697}
\email{bozbagci@math.uci.edu}
%\subjclass{57R55, 57R65, 57M12}
\date{\today}

\begin{abstract}
Let $M$ be a smooth 4-manifold which admits a
genus $g$ $C^{\infty}$-Lefschetz fibration over $D^2$ or $S^2$. 
We develop a technique to
compute the signature of $M$ using 
the global monodromy of this fibration.
As a corollary we prove that there is no hyperelliptic 
Lefschetz fibration (of any genus) over $S^2$ 
with only reducible singular fibers. 
We also prove that for each $g \geq 1$, there exist 
$$k_g \leq 10- \frac{6g+4}{g^2} \; \; 
and \; \; l_g \leq 2g-10+ \frac{4g+4}{g^2} $$ 
such that
if a $4$-manifold admits a 
hyperelliptic $C^{\infty}$-Lefschetz
fibration of genus $g$ over $S^2$ 
then its signature $\sigma$ and 
its Euler characteristic $e$ 
satisfy the inequality $c_1^2 \leq k_g \chi + l_g $, 
where $c_1^2 =3\sigma+2e$ and $\chi=\frac{1}{4}(\sigma+e)$. 
In particular, we show that we can choose $k_2=6$, $ l_2=-3 $,
$k_3=7.25$, $l_3=-2.75$, $k_4=8.25$ and $l_4=-0.75$. 

\end{abstract}

\maketitle
\setcounter{section}{-1}
\section{Introduction}

The signature of a smooth 4-manifold which admits a  hyperelliptic 
Lefschetz fibration of genus $g$
over a closed surface can be computed using the 
{\em local signature formula}, given by Matsumoto ({\cite{m1}}, 
{\cite{m2}}) for
$g=1,2$ and more recently extended by Endo {\cite{e}} for $g \geq 3$.

In this paper we present a method to compute 
the signature of a smooth 4-manifold which 
admits an arbitrary (not necessarily hyperelliptic) 
Lefschetz fibration of any genus
over $D^2$ or $S^2$.   
A Lefschetz fibration on a smooth 4-manifold
$M$ gives rise to a handlebody description of $M$, which is 
determined by a sequence of {\em vanishing cycles}. 
We use this handlebody 
description {\cite{k}} and 
{\em Wall's nonadditivity formula} for signatures {\cite{w}}
in order to compute the signature of $M$. 
Hence we calculate a {\em `relative signature'} 
corresponding to each singular fiber of the given fibration on $M$.

Despite the fact that the vanishing cycles are defined
up to isotopy, our technique shows that the signature of
a 4-manifold which admits a Lefschetz fibration 
depends only on the algebraic data given 
by the homology classes of the vanishing cycles.

Recent results in symplectic topology shows that 
Lefschetz fibrations provide a topological
characterization of symplectic 4-manifolds: Donaldson \cite{d}
has shown that, after perhaps blowing up, a closed symplectic\;
4-manifold admits a Lefschetz fibration over $S^2$, and conversely Gompf 
\cite{gs} has shown that most Lefschetz fibrations 
are symplectic -- the exceptions all have fiber-genus 1 and are 
blow-ups of torus fibrations with no critical points. 
Hence by computing the signatures of Lefschetz fibrations
of any genus we hope to answer some of the problems
in the {\em geography} of symplectic $4$-manifolds 
({\cite{s1}}, {\cite{s2}}, {\cite{gs}}).

It is well-known that the complex surfaces satisfy the 
BMY-inequality $c_1^2 \leq 9\chi$. In {\cite{s2}}, Stipsicz 
proves that $0 \leq c_1^2 \leq 10 \chi$ for
the (relatively minimal) genus $g (\geq 1)$ Lefschetz fibrations over
closed surfaces of nonzero genus. His result, however, does not extend to 
cover the fibrations over $S^2$. Hence the symplectic 
version of the BMY-inequality
is still missing. Stipsicz {\cite{s2}} also points out that 
$c_1^2 \leq 10 \chi +2g-1$ holds for any Lefschetz fibration.

Lefschetz fibrations on smooth 4-manifolds are generalized from
holomorphic Lefschetz fibrations on complex surfaces. It is 
conjectured by G.Tian that any irreducible hyperelliptic (smooth) Lefschetz 
fibration is holomorphic.
Here irreducible means that 
the global monodromy of the fibration does not 
contain a Dehn twist about a seperating curve. 
There are, however, noncomplex smooth 4-manifolds which admit 
reducible hyperelliptic Lefschetz fibrations (\cite{fs}, \cite{os}).

We note that the signature of an irreducible hyperelliptic Lefschetz 
fibration is always negative. 
   
{\it Acknowledgement.} \ The author would like to thank John Etnyre, 
Terry Fuller, Ludmil Katzarkov, Ron Stern and Richard Wentworth 
for helpful conversations. In particular, the author 
would like to thank Terry Fuller for stimulating his 
interest in Lefschetz fibrations. The author would also
like to thank Yukio Matsumoto for informing him about the preprint
{\cite{e}} by Hisaaki Endo. The author would like to thank 
Andras Stipsicz for commenting on earlier versions of this paper.
Finally, the author wishes to express his gratitute to Ron Stern 
for encouragement and many useful discussions.

\section{Preliminaries}

\subsection{Mapping Class Groups}

{\Def Let $\S_g$ be a closed oriented surface of genus $g$. 
Let ${Diff}^+ (\S_g)$ be the group of all orientation 
preserving self diffeomorphisms of $\S_g$.
Let ${Diff}_0^+ (\S_g)$ be the subgroup of ${Diff}^+ (\S_g)$    
consisting of all self diffeomorphisms isotopic to the identity.
Then we define the {\em mapping class group} of genus $g$ as
$$ {\CM}_g = {Diff}^+ (\S_g) /{Diff}_0^+ (\S_g).$$ 

The {\em hyperelliptic mapping class group} 
${\CH}_g $ of genus $g$ is defined ({\cite{e}}) as
the subgroup of $ {\CM}_g $ which consists of all isotopy classes 
commuting with the isotopy class of the hyperelliptic involution 
$\iota :\Sigma _g\longrightarrow
\Sigma _g$  .}
 
It is known ({\cite{bh}}) that 
the hyperelliptic mapping class group ${\CH}_g $ 
agrees with the mapping class group ${\CM}_g $ for $g=1,2$.
 
We will denote a positive Dehn twist about 
a simple closed curve $\g$
by $D(\g) \in {\CM}_g $ and 
we will use the functional notation for the products in
${\CM}_g$, e.g., $D(\b)D(\a)$ will denote the composition where
we apply $D(\a)$ first and then $D(\b)$.

\subsection{Smooth Lefschetz Fibrations}

{\Def Let $M$ be a compact, oriented smooth 4-manifold, and
let $B$ be a compact, oriented 2-manifold. A proper smooth
map $f:M\to B$ is a {\em ($C^{\infty}$-) Lefschetz fibration}
if there exist points $b_1,\ldots,b_m \in {\mathrm interior}
(B)$ such that
\begin{itemize}
\item[(1)] $\{b_1,\ldots,b_m \}$ are the critical values of
$f$, with $p_i\in f^{-1}(b_i)$ a unique critical point
of $f$, for each $i$, and
\item[(2)] about each $b_i$ and $p_i$, there are complex 
coordinate neighborhoods such that locally $f$ can be expressed
as $f(z_1,z_2)=z_1^2+z_2^2$.
\end{itemize}}

It is a consequence of this definition that
$$f|_{f^{-1}(B-\{b_1,\ldots,b_m \})}:
f^{-1}(B-\{b_1,\ldots,b_m \})\to B-\{b_1,\ldots,b_m \}$$
is a smooth fiber bundle over $B-\{b_1,\ldots,b_m \}$
with fiber diffeomorphic to a 2-manifold $\S_g$, and so
we refer to $f$ (and sometimes also the manifold $M$) 
as a {\em genus $g$ Lefschetz fibration} (or
a {\em Lefschetz fibration of genus $g$}).  
Two genus $g$ Lefschetz fibrations $f:M\to B$ and
$f^{\prime}:M^{\prime}\to B^{\prime}$ are {\em equivalent}
if there are diffeomorphisms $\Phi:M\to M^{\prime}$ and
$\phi:B\to B^{\prime}$ such that 
$f^{\prime}\Phi=\phi f.$

We always assume that our Lefschetz fibrations are 
{\em relatively minimal}, namely that no fiber contains
an embedded 2-sphere of self-intersection number $-1$.
We also assume that there is at least one singular fiber 
in each fibration.

If $f:M\to D^2$ is a smooth genus $g$ Lefschetz fibration,
then we can use the Lefschetz fibration to produce a handlebody
description of $M$. 
We select a regular value $b_0 \in {\mathrm interior}(D^2)$
of $f$, an identification $f^{-1}(b_0)\cong \S_g$,
and a collection of arcs $s_i$ in
${\mathrm interior}(D^2)$ with each $s_i$
connecting $b_0$ to $b_i$, and otherwise disjoint
from the other arcs. We also assume that the critical
values are indexed so that the arcs $s_1,\ldots,s_m$
appear in order as we travel counterclockwise in a 
small circle about $b_0$.
Let $V_0,\ldots,V_m$ denote a collection of small disjoint
open disks with $b_i \in V_i$ for each $i$.

To build our description of $M$, we observe first
that $f^{-1}(V_0)\cong \S_g\times D^2$,
with $\bdy V_0\cong \S_g\times S^1.$
Enlarging $V_0$ to include the critical value
$b_1$, it can be shown that
$f^{-1}(V_0\cup \nu(s_1) \cup V_1)$
is diffeomorphic to $\S_g\times D^2$
with a 2-handle $h_1$ attached along a circle
${\g_1}$ contained in a fiber
$\S_g\times pt. \subset \S_g\times S^1.$
Moreover, condition (2) in the definition of a
Lefschetz fibration requires that $h_1$
is attached with a framing $-1$ relative
to the natural framing on ${\g_1}$
inherited from the product structure of $\bdy V_0$.
${\g_1}$ is called a vanishing cycle.
In addition, $\bdy((\S_g\times D^2)\cup h_1)$ is
diffeomorphic to a $\S_g$-bundle over $S^1$ whose monodromy
is given by $D({\g_1})$, a positive Dehn twist about
${\g_1}$. Continuing counterclockwise about $b_0$,
we add the remaining critical values to our description, 
yielding that
$$M_0\cong f^{-1}(
V_0\cup (\bigcup_{i=1}^m \nu(s_i))
\cup (\bigcup_{i=1}^m V_i))$$
is diffeomorphic to $(\S_g\times D^2)\cup 
(\bigcup_{i=1}^m h_i)$, where each $h_i$
is a 2-handle attached along a vanishing cycle ${\g_i}$
in a $\S_g$-fiber in $\S_g\times S^1$ with relative framing $-1$.
Furthermore, 
$$\bdy M_0 \cong
\bdy((\S_g\times D^2)\cup (\bigcup_{i=1}^m h_i))$$
is a $\S_g$-bundle over $S^1$ with monodromy given
by the composition $D({\g_m})\cdots D({\g_1})$.
We will refer to the cyclically ordered collection 
$(D({\g_1}) ,..., D({\g_m}))$ (or the product 
$D({\g_m})\cdots D({\g_1})$)
as the {\em global monodromy} of this fibration. 

We can extend this description to Lefschetz fibrations over $S^2$
as follows:

Assume that $f:M\to S^2$ is a smooth genus $g$ 
Lefschetz fibration.
Let $M_0=M-\nu(f^{-1}(b))$, where
$\nu(f^{-1}(b))\cong \S_g\times D^2$ denotes a regular neighborhood 
of a nonsingular fiber $f^{-1}(b)$.
Then $f|_{M_0}:M_0 \to D^2$ is a smooth Lefschetz fibration.
If $(D({\g_1}) ,..., D({\g_m}))$ is the global monodromy
of the fibration $f|_{M_0}:M_0 \to D^2$, 
then $D({\g_m})\cdots D({\g_1})$
is isotopic to the identity since also $\bdy M_0 \cong \S_g\times S^1$. 
Finally, to extend our description
of $M_0$ to $M$, we reattach $\S_g\times D^2$ to 
$(\S_g\times D^2)\cup (\bigcup_{i=1}^m h_i)$
via a $\S_g$-fiber preserving map of the boundary.
This extension is unique up to equivalence
for $g\geq 2$ \cite{k}.

{\Rem Although the description of the monodromy
corresponding to each individual critical value $b_i$
as a Dehn twist depends on the choice of arc $s_i$,
other choices of arcs (and of the central identification
$f^{-1}(x_0)\cong \S_g$) do not change the Lefschetz
fibration on $M_0$, up to equivalence {\cite{gs}}}. 

{\Def Let $f:M\to S^2$ be a smooth genus $g$
Lefschetz fibration with global monodromy 
$(D({\g_1}) ,..., D({\g_m}))$. We will call $f:M\to S^2$ 
a {\em hyperelliptic Lefschetz fibration} of genus $g$ iff there 
exists $h \in {\CM}_g$ such that $h D({\g_i}) h^{-1} \in 
{\CH}_g$ for all i, $1 \leq i \leq m$.}

{\Rem All Lefschetz fibrations of genus one and genus two 
are hyperelliptic since 
${\CH}_g = {\CM}_g $ for $g=1,2$. }

\subsection{Wall's Non-Additivity Formula}

\indent

If two compact oriented $4$-manifolds are glued
by an orientation reversing diffeomorphism
of their boundaries, then the signature of their 
union is the sum of their signatures. This is known
as the Novikov additivity. But it is often
desirable to consider the more general case of
gluing: along a common submanifold, which may 
itself have boundary, of the boundaries of the
original manifolds. However, the Novikov
additivity does not hold in this general case.
Wall {\cite{w}} derives a formula for the deviation
from additivity in the general case, which is known 
as the Wall's nonadditivity formula. 

We will give a specific case of his formula:

Let $X_-$, $X_0$, $X_+$ be $3$-manifolds and $Y_-$ and
$Y_+$ be $4$-manifolds such that

$\bdy X_-=\bdy X_0=\bdy X_+=Z$, 

$\bdy Y_-=X_- \cup X_0 $, 

$\bdy Y_+= X_0 \cup X_+ $;

write $Y=Y_- \cup Y_+$ and $X=X_- \cup X_0 \cup X_+$. (Figure $1$)

\begin{figure}[h]
 \begin{center}
    \includegraphics{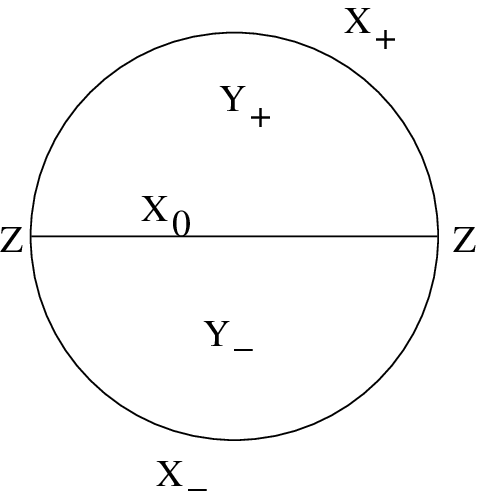}
  \caption{} \label{wall}
  \end{center}
 \end{figure}

Suppose that $Y$ is oriented inducing orientations of
$Y_-$ and $Y_+$. Orient the rest so that 

$\bdy_* [Y_-]=[X_0]-[X_-]$, 

$\bdy_* [Y_+]=[X_+]-[X_0]$,

$\bdy_* [X_-]=\bdy_* [X_0]=\bdy_* [X_+]=[Z]$.

Write $V=H_1 (Z;\IR)$; let $A,B,C$ be the kernels 
of the maps induced by the inclusions of $Z$ in 
$X_-$,$X_0$ and $X_+$ respectively. Then
$dimA=dimB=dimC= \frac{(dim V)}{2}$.

Let $\Phi$ denote 
the oriented intersection numbers in $Z$. 
Let $W= \frac{C \cap (A+B)}{(C \cap A)+(C \cap B)}$.
Then Wall defines a symmetric bilinear map
$\Psi :W \times W \rightarrow \IR$ as follows:
First define $\Psi^\prime :C \cap (A+B) \times
C \cap (A+B) \rightarrow \IR$ by 
$$\Psi^\prime (c,c^\prime)=\Phi(c,a^\prime)$$
where $a^\prime+b^\prime+c^\prime=0$ for some $b^\prime \in B$.
Then $\Psi^\prime$
induces a map $\Psi$ on $W$.

The signature of the symmetric bilinear map $\Psi$
will be denoted by $\sigma(V;C,A,B)$.

We also denote the signature of a 4-manifold M
as $\sigma(M)$ in the rest of this paper.

We are now ready to state Wall's formula:

{\Thm {\cite{w}}
$\sigma (Y)=\sigma (Y_-)+ \sigma (Y_+)-\sigma(V;C,A,B)$.}

\subsection{Local Signature Formula}

\indent 

The following theorem was proven by Matsumoto for $g=1,2$ using 
the fact that the cohomology 
class of Meyer's signature cocyle has finite order in the 
cohomology group $H^2 ({\CM}_g , \IZ)$. Recently, Endo proved 
the $g \geq 3$
case by observing the finiteness of the order of the 
cohomology class of the signature cocycle restricted to the
hyperelliptic mapping class group ${\CH}_g$.   

{\Thm \label{frac}
$(${\cite{m1}}, {\cite{m2}}, 
{\cite{e}}$)$ Let $M$ be a 
$4$-manifold which admits a 
hyperelliptic Lefschetz fibration of genus g over $S^2$.
Let $n$ and $s=\sum_{h=1}^{[\frac{g}{2}]}s_h$
be the numbers of nonseparating and separating
vanishing cycles in the global monodromy of this 
fibration, respectively. Then
$$\sigma(M)=-\frac{g+1}{2g+1}n+
\sum_{h=1}^{[\frac{g}{2}]}(\frac{4h(g-h)}{2g+1} -1)s_h .$$}

{\Rem Here $s_h$ denotes the number of separating vanishing cycles 
which separate the genus $g$ surface into two surfaces one of which 
has genus $h$.} 

\section{Main Theorems}
\indent
In this section we explain our main idea and establish the main 
theorems to develop an algorithm to compute the signature of a 
a 4-manifold which admits a Lefschetz fibration over $D^2$ or $S^2$ 
using the global monodromy of this fibration. 

{\Def Let $X$ be a $4$-manifold with boundary $\bdy X\cong
\S_g\times I/(x,1)\sim(\phi(x),0)$, where $\phi$ 
is a self-diffeomorphism of $\S_g$. Let $X^\prime$ denote
the resulting $4$-manifold 
after attaching a $2$-handle to $X$ 
along a simple closed curve 
${\g}$ on $\S_g\times\{pt\}$ with framing $-1$ (relative
to the product framing).  
Then $\sigma(\phi,\g)$ is defined as $\sigma(X^\prime) -
\sigma(X)$. }

{\Thm \label{burak1}
Let $M$ be a $4$-manifold which admits a genus $g$ Lefschetz 
fibration over $D^2$ or $S^2$. 
Let $(D({\g_1}),..., D({\g_t}))$ be
the global monodromy of this fibration. Let $ D({\g_0})$
denote the identity map. Then 
$$\sigma(M)=\S_{i=1}^t \sigma(D(\g_{i-1})\cdots D(\g_0), \g_i ),$$
where $ \sigma(D(\g_{i-1})\cdots D(\g_1), \g_i ) 
\in \{-1,0,+1\}$ for all i, $1\leq i \leq t$.}

\begin{proof}
It suffices to prove the result for Lefschetz 
fibrations over $D^2$. (By Novikov additivity it extends to
Lefschetz 
fibrations over $S^2$.)  
We use the handlebody description of $M$
and Wall's formula as follows:

We start with a copy of $M_0=\S_g\times D^2$. We attach 
a $2$-handle to $M_0$ along $\g_1$ with framing $-1$.
Let $M_1$ denote the resulting manifold. Then $\bdy M_1$ will 
have monodromy $D({\g_1})$, a positive Dehn twist about $\g_1$.
Now we attach
another $2$-handle to $M_1$ along $\g_2$. Let $M_2$
denote the resulting manifold. Proceeding in this manner
we get the manifolds  $M_1,M_2,..., M_t$. 

We are going to apply Wall's formula at each step
of this contruction to compute the signature of $M$.
In order to apply Wall's formula we set up
the following notation:

$\phi$: a self-diffeomorphism of $\S_g$  

$X$: a $4$-manifold with boundary $\bdy X\cong 
\S_g\times I/(x,1)\sim(\phi(x),0)$.

$\g$: a simple closed curve embedded in a fiber $\S_g\times\{pt\}$. 

$X^\prime$: resulting $4$-manifold
after attaching a $2$-handle to $X$
along a simple closed curve
${\g}$ on $\S_g\times\{pt\}$ with framing $-1$ (relative
to the product framing). 

$\nu(\g)$: a regular neighborhood of $\g$ in $\bdy X$.

$i_*$: the induced map on the homology 
by the inclusion of appropriate spaces.

Now we define $Y_+,\;Y_-,\;X_+,\;X_0,\;X_-,\;Z$ in Wall's formula as follows:

$Y_-=D^2 \times D^2$, $\;Y_+=X$,

$\bdy Y_-=\bdy{(D^2 \times D^2)}=S^1 \times D^2 \cup D^2 \times S^1$,

$\bdy Y_+=\S_g\times I/(x,1)\sim(\phi(x),0)$,

$X_0=S^1 \times D^2 \cong \nu(\g)$,
$\;X_-=D^2 \times S^1 $, $\;X_+=\bdy X-\vo(\g)$,

$Z=S^1 \times S^1 \cong \bdy\nu(\g) \cong \bdy{(\bdy X - \vo(\g))}.$

Hence,

$A=Ker(i_*: H_1(S^1\times S^1;{\IR}) \rightarrow H_1(D^2\times S^1;{\IR})),$

$B=Ker(i_*: H_1(S^1\times S^1;{\IR}) \rightarrow H_1(S^1\times D^2;{\IR})),$

$C=Ker(i_*: H_1(\bdy \nu(\g);{\IR}) 
\rightarrow H_1(\bdy X - \vo(\g);{\IR})).$

Let $l$ be the longitude $S^1 \times \{pt\}$ and $m$ be 
the meridian $\{pt\} \times \bdy D^2$ of $X_0=S^1 \times D^2$. 
Then $A=<[l]>$ and $B=<[m]>$. We also know 
that $C$ is a $1$-
dimensional subspace of $$H_1(S^1 \times S^1;{\IR})=<[l],[m]> \cong 
{\IR}^2.$$

Let $\Phi$ be the intersection form on $Z=S^1 \times S^1$ 
and $W= \frac{C \cap (A+B)}{(C \cap A)+(C \cap B)} = 
\frac{C}{(C \cap A)+(C \cap B)}$. Hence
$W=\{0\}$ if $C=A$ or $C=B$ and $W=C$ otherwise.
Now assume that $C\neq A$ and $C\neq B$. Then    
$C=<c>=<p[l]+q[m]>$ for some $p,q \in \IR$ and
$\Psi(c,c)=\Phi(c,a^\prime)$ where
$c+a^\prime +b^\prime=0$ for some $a^\prime \in A$ and
$b^\prime \in B$. ($\Psi$ is the bilinear form in Wall's formula). 
Let $a^\prime=-p[l]$ and $b^\prime=-q[m]$. 
Then we have, 
$$\Psi(c,c)=\Phi(c,-p[l])=\Phi(p[l]+q[m],-p[l])=-pq \Phi([m],[l])=pq.$$
Therefore signature of $\Psi$ is given by the sign of $pq$.
  
Hence by Wall's formula 

$\sigma(X^\prime)=\sigma(X)+\sigma(D^2 \times D^2)
-\sigma({\IR}^2;C,A,B)$

$=\sigma(X)-signature(\Psi)=\sigma(X)-sign(pq)$.

This proves the theorem by setting $X= M_i$ for $i=1,2,...,t-1$.

\end{proof}

So the idea to compute the signature of a genus $g$ 
Lefschetz fibration
is very simple. For each 2-handle that we attach to 
$\S_g\times D^2$ along a vanishing cycle, there is a corresponding
{\em relative signature} $\in \{-1,0,+1\}$. Once we attach all
the 2-handles, the sum of the relative signatures will be signature 
of the 4-manifold. The difficulty is to compute the relative signatures
using the vanishing cycles (or more precisely using only the homology 
classes of the vanishing cycles). The following technical 
theorems will be helpful in 
computations.

{\Thm \label{local}
In addition to the notation above, 
let $\{a_1,b_1,a_2,b_2,...,a_g,b_g\}$ be the standard basis
for $H_1(\S_g;\IR)$. (We will 
use the letters $a_i$ and $b_i$ 
also to denote the curves which represent
the homology classes $a_i$ and $b_i$, 
respectively, for $1 \leq i \leq g$.) Then

\begin{itemize}

\item[(1)]If $\g$ is a nonseparating curve, then there
exists a longitude $l^\prime$ and a meridian $m^\prime$ 
of $\bdy (\bdy X -\vo(\g))$
such that 
$$i_*[l^\prime]=[\g]\in H_1(\bdy X - \vo(\g);{\IR})$$
$$i_*[m^\prime]= \frac{e -\phi_*(e)}{e.[\g]} 
\in H_1(\bdy X - \vo(\g);{\IR})$$
for all $e \in\{a_1,b_1,a_2,b_2,...,a_g,b_g\}$,
where $e.[\g] \neq 0$.
\item[(2)]If $\g$ is a separating curve, then 
$\sigma(X^\prime)=\sigma(X)-1$ , i.e., $\sigma(\phi,\g) =-1$. 

\end{itemize}}

\begin{proof}

We recall that $\bdy X$ is a mapping torus, i.e.,
$\bdy X \cong \S_g\times I/(x,1)\sim(\phi(x),0)$ and 
$\g$ is a curve on a fiber $\S_g \times \{pt\}$.
We note that a regular neighborhood of $\g$ in $\S_g$ 
is given by $\g \times I_1$. Hence a regular 
neighborhood of $\g$ in $\bdy X$ 
is given by $\g \times I_1 \times I_2$ where
$I_2$ is a small neighborhood of the $\{pt\}$
in $S^1 = I/(1 \sim 0)$. This neighborhood 
of $\g$ is called the product 
neighborhood {\cite{k}}. 

Now let us push off $\g$ to the boundary of $\bdy X - \vo(\g)$.
Denote the push off of $\g$ as $l^\prime$. Moreover 
if we identify $I_1\times I_2$ as $D^2$ and denote
$\bdy D^2$ as $m^\prime$, then $\{ l^\prime, m^\prime\}$
will be a longitude-meridian pair for $\bdy(\bdy X - \vo(\g))$.
Then clearly  
$$i_*[l^\prime]=[\g]\in H_1(\bdy X - \vo(\g);{\IR})$$ 

On the other hand, to find the image of $m^\prime$
we observe the following:

Assume that $e.[\g] =1$ for 
some $e\in \{a_1,b_1,a_2,b_2,...,a_g,b_g\}$.
Then we locally have the picture in Figure 2 in a 
neighborhood of the point where $e$ and $\g$ meet.

\begin{figure}[h]
 \begin{center}
    \includegraphics{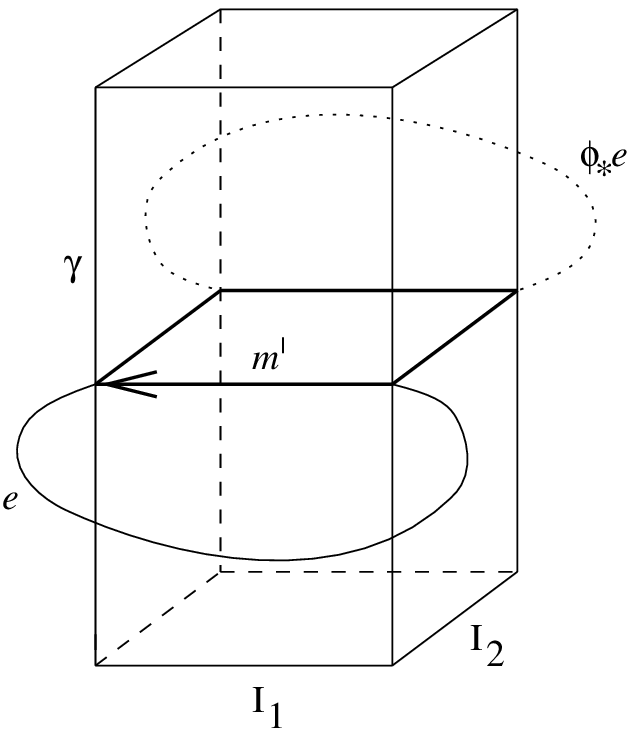}
  \caption{} \label{local1}
  \end{center}
 \end{figure}

This proves that  
$$i_*[m^\prime]=e -\phi_*(e)\in H_1(\bdy X - \vo(\g);{\IR}).$$

Now assume that $e.[\g] =-1$ 
for some $e\in \{a_1,b_1,a_2,b_2,...,a_g,b_g\}$.
Then we locally have a similar  picture in a 
neighborhood of the point where $e$ and $\g$ meet, except for the orientations.

This proves that  
$$i_*[m^\prime]=\phi_*(e) - e \in H_1(\bdy X - \vo(\g);{\IR}).$$

Since these are local results it follows 
combining these two observations
that
$$i_*[m^\prime]= \frac{e -\phi_*(e)}{e.[\g]} 
\in H_1(\bdy X - \vo(\g);{\IR}).$$

To prove the second part of the theorem we note that if 
$\g$ is a separating curve in $\S_g$ then it is homologically
trivial. Thus $i_*[l^\prime]=0$. This implies 
that $Ker(i_*)=<[l^\prime]>$. 

Note that, in terms of the 
bases $\{[l],[m]\}$ of $H_1(\bdy(S^1 \times D^2);{\IR})=
H_1(S^1 \times S^1;{\IR})$ and 
$\{[l^\prime] ,[m^\prime] \}$ of $H_1(\bdy(\bdy X - \vo(\g));{\IR})$, 
attaching a $2$-handle by $-1$ framing means that 
we identify $[l]$ with $[l^\prime]-[m^\prime]$ and $[m]$ with $[m^\prime]$.

So if we transform the $Ker(i_*)$ to the $\{[l],[m]\}$ plane we see that
$Ker(i_*)=C=W=<[l]+[m]>$ which implies that $\sigma(X^\prime)=
\sigma(X)-(+1)$ (cf. Theorem~\ref{burak1}).

\end{proof}

{\Prop \label{homology}
We use the same notation as in Theorem~\ref{local}. 

\begin{itemize}

\item[(1)] Let $\g =a_i$ for some i, $1\leq i \leq g$. 
Then $H_1(\bdy X - \vo(a_i);{\IR})=$ 

$\langle a_1,b_1,a_2,b_2,...,a_g,b_g,
b_i^\prime,t$ $|$ $a_j=\phi_* a_j $ for 
all j, $b_j=\phi_* b_j$ for
all $j \neq i,$ $b_i^\prime=\phi_* b_i  \rangle$. 

Moreover 
$i_*[l^\prime]=a_i $ and $i_*[m^\prime]=b_i - b_i^\prime $.

\item[(2)] Let $\g =b_i$ for some i, $1\leq i \leq g$. Then
$H_1(\bdy X - \vo(b_i);{\IR})=$

$ \langle a_1,b_1,a_2,b_2,...,a_g,b_g,
a_i^\prime,t$ $|$ $b_j=\phi_* b_j $ for
all j, $a_j=\phi_* a_j$ for
all $j \neq i,$ $a_i^\prime=\phi_* a_i  \rangle$. 

Moreover
$i_*[l^\prime]=b_i $ and $i_*[m^\prime]=a_i^\prime - a_i  $.

\end{itemize}}

\begin{proof}
Assume that $\g =a_i$ for some i, $1\leq i \leq g$. 
We first use Van-Kampen's theorem to compute ${\pi}_1(\bdy X - \vo(a_i))$. 
Write $\bdy X - \vo(a_i)=E_1 \cup E_2$ as follows:
Let $E_1=\S_g\times[0,1/2]$ and $E_2=\S_g\times[1/2,1]$.  
Then glue $\S_g\times\{1/2\} \subset E_1$ with
$\S_g\times\{1/2\} \subset E_2$ by the identity map
except a neighborhood of $a_i$, namely $a_i \times I \subset \S_g$.
Denote the result as $E^\prime$. 
By a trivial calculation we get the following 
presentation: 

$\pi_1(E^\prime)=
\langle a_1,b_1,a_2,b_2,...,a_g,b_g,
b_i^\prime$ $|$ $\prod_{j=1}^g$ 
$[a_j,b_j],$ $[a_i,b_i^\prime]
\prod_{j\neq i} [a_j,b_j] \rangle$

Finally we abelianize this presentation 
after gluing  $\S_g\times\{0\} \subset E_1$ with
$\S_g\times\{1\} \subset E_2$ using the map $\phi$
to get $\bdy X - \vo(a_i)$.

$i_*[l^\prime]=a_i $ and $i_*[m^\prime]=b_i - b_i^\prime $
follows from Theorem~\ref{local} because $a_i$ intersects $b_j$
only once iff $i=j$.

Second part is obtained similarly.

\end{proof}

\section{Applications}

First we give an immediate application of the theorems 3 and 4. 

{\Cor \label{bound}
Let $M$ be a $4$-manifold which admits
a genus g Lefschetz fibration over $D^2$ or $S^2$.
Let $n$ and $s$ be the numbers of nonseparating and 
separating vanishing cycles in the global monodromy of this 
fibration, respectively. Then $\sigma(M) \leq n-s .$}

\begin{proof}

Suppose that we build up the 4-manifold $M$ from $\S_2 \times D^2$
by attaching 2-handles. By Theorem~\ref{local}, every time we attach a
$2$-handle along a separating curve the signature
of the resulting 4-manifold will be one less than
the signature of the 4-manifold before we attach the 2-handle.
Thus Theorem~\ref{burak1} implies the 
upper bound $ n-s$ on the signature.

\end{proof}

{\Cor \label{reducible}
There is no hyperelliptic Lefschetz fibration (of any genus) over $S^2$ 
with only reducible singular fibers. (Here reducible means 
that the local monodromy
corresponding to the singular fiber is a Dehn twist about 
a separating curve.)}

\begin{proof}

Let $M$ be a $4$-manifold which admits a Lefschetz
fibration of genus $g$ over $S^2$
with global monodromy $(D({\g_1}),...,D({\g_s}))$ , where
$s=\sum_{h=1}^{[\frac{g}{2}]}s_h$ and
${\g_i}$ is separating for each i, $ 1 \leq  i \leq s $.
Then, by the local signature formula,
\[ \sigma(M) =\left\{ \begin{array}{ll}
\sum_{h=1}^{[\frac{g}{2}]} (\frac{4h(g-h)}{2g+1} -1) s_h \geq 0 & 
\mbox{if $g \geq 3$} \\
-s/5 & \mbox{if $g=2$} 
\end{array}
\right. \]

But on the other hand
$\sigma(M)=-s$ according
to Theorem~\ref{local}.
Hence $s=0$. (This is trivially true for $g=1$ since any vanishing 
cycle on a torus is
nonseparating.) This proves the desired result
since we assume (by definition) that there exists at least one
singular fiber in each Lefschetz fibration.

\end{proof}

Next we combine our results with the local signature 
formula for the hyperelliptic Lefschetz 
fibrations to give an upper bound for the signatures
of these fibrations.

{\Cor \label{hbound}
Let $M$ be a $4$-manifold which admits
a hyperelliptic Lefschetz fibration of genus $g$
over $S^2$.
Let $n$ and $s$
be the numbers of nonseparating and separating
vanishing cycles in the global monodromy of this
fibration, respectively. Then $\sigma(M) \leq n-s-4 .$}

\begin{proof}

We first note that we can improve the inequality
$$\sigma(M) \leq n-s$$
given in Corollary~\ref{bound} to
$$\sigma(M) \leq n-s-1$$
for hyperelliptic Lefschetz fibrations as follows:

Suppose that we attach the first
2-handle along a nonseparating curve. 
We can always assume this because $n \geq 1$ (since 
we proved in Corollary~\ref{reducible} that $n \neq 0$) and we 
can cyclically permute 
the vanishing cycles in the global monodromy of a Lefschetz fibration. 
Moreover we can easily show that
if we start attaching handles
along a nonseparating curve
then the signature of the resulting 4-manifold (after
attaching the very first handle)
will be the same as $\sigma(\S_2 \times D^2)$, which is zero.

Next note that $\sigma(M) \leq n-s-1$ is equivalent to
$$4\sum_{h=1}^{[\frac{g}{2}]}h(g-h)s_h \leq (3g+2)n-(2g+1)$$
using the local signature formula.

Assume that $g$ is odd. Endo {\cite{e}} proves that 
$$n+4 \sum_{h=1}^{[\frac{g}{2}]}h(2h+1)s_h  \equiv 0  \;\; (\bmod 4(2g+1)) .$$
Hence $$n=4c(2g+1)-4 \sum_{h=1}^{[\frac{g}{2}]}h(2h+1)s_h$$
for some integer c. Substituting into the inequality above 
(and dividing by 4) we get
$$\sum_{h=1}^{[\frac{g}{2}]}h(g-h)s_h \leq 
(3g+2)[c(2g+1)- \sum_{h=1}^{[\frac{g}{2}]}h(2h+1)s_h]-\frac{1}{4}(2g+1).$$
Hence
$$\sum_{h=1}^{[\frac{g}{2}]}h(g-h)s_h \leq
(3g+2)[c(2g+1)- \sum_{h=1}^{[\frac{g}{2}]}h(2h+1)s_h]-\frac{1}{4}(2g+2)$$
since $2g+1 \equiv 3 \;\;  (\bmod 4)$ .

But this inequality, in turn, implies that 
$$4\sum_{h=1}^{[\frac{g}{2}]}h(g-h)s_h \leq (3g+2)n-(2g+2)$$ 
which is equivalent to
$$\sigma(M) \leq n-s-1-\frac{1}{2g+1}.$$
Since $\sigma(M)$ is an integer,
$$\sigma(M) \leq n-s-2.$$
Iterating the same argument, we obtain 
$$\sigma(M) \leq n-s-4.$$
(We use $2(2g+1) \equiv 2 \;\;  (\bmod 4)$ and 
$3(2g+1) \equiv 1 \;\;  (\bmod 4))$.

Similarly, if $g$ is even, then one can  
use the corresponding result by Endo:
$$n+4 \sum_{h=1}^{[\frac{g}{2}]}h(2h+1)s_h  \equiv 0 \;\;  (\bmod 2(2g+1)) .$$
(Note that $2(3g+2)(2g+1) \equiv 0 \;\;  (\bmod 4)$, if g is even.)
\end{proof}

The following corollary is our main result concerning 
the {\em geography} of the hyperelliptic Lefschetz
fibrations.

{\Cor \label{gengeo}
For each $g \geq 1$, there exist 
$$k_g \leq 10- \frac{6g+4}{g^2} \; \; 
and \; \; l_g \leq 2g-10+ \frac{4g+4}{g^2} $$ 
such that
if a $4$-manifold admits a 
hyperelliptic Lefschetz
fibration of genus $g$ over $S^2$ 
then its signature $\sigma$ and 
its Euler characteristic $e$ 
satisfy the inequality $c_1^2 \leq k_g \chi + l_g $, 
where $c_1^2 =3\sigma+2e$ and $\chi=\frac{1}{4}(\sigma+e)$. }

\begin{proof}

Let $M$ be a $4$-manifold which admits a
hyperelliptic Lefschetz fibration of genus $g$ over $S^2$.
Let $n$ and $s=\sum_{h=1}^{[\frac{g}{2}]}s_h$ 
be the numbers of nonseparating and separating
vanishing cycles in the global monodromy of this 
fibration, respectively. Then 
$$\sigma(M)=-\frac{g+1}{2g+1}n+
\sum_{h=1}^{[\frac{g}{2}]}(\frac{4h(g-h)}{2g+1} -1)s_h $$ 
and 
$$e(M)=n+\sum_{h=1}^{[\frac{g}{2}]}s_h - 4(g-1).$$
Thus 
$$ c_1^2 (M) = 3\sigma(M) + 2e(M) = \frac{1}{2g+1} [(g-1)n-8(g-1)(2g+1)+
\sum_{h=1}^{[\frac{g}{2}]}(12h(g-h)-2g-1)s_h  ]$$
and
$$ \chi (M) =\frac{\sigma(M) +e(M)}{4} = \frac{1}{4(2g+1)}[gn-4(g-1)(2g+1)+
\sum_{h=1}^{[\frac{g}{2}]}4h(g-h)s_h ].$$
We want to find $k_g$ and $l_g$ such that 
$c_1^2 (M) \leq k_g \chi (M) + l_g .$
Using the equalities above we get,
$$\sum_{h=1}^{[\frac{g}{2}]}[(12-k_g)h(g-h)-2g-1]s_h  \leq 
[\frac{k_g}{4}g -(g-1)]n+ m_g ,$$  
where $m_g = -(g-1)(2g+1)k_g +(2g+1)l_g +8(g-1)(2g+1).$
Now we note that 
$$\sum_{h=1}^{[\frac{g}{2}]}[(12-k_g)h(g-h)-2g-1]s_h =
\sum_{h=1}^{[\frac{g}{2}]}h(g-h)(12- k_g - \frac{2g+1}{h(g-h)})s_h. $$
If $g$ is even, then 
$$\sum_{h=1}^{[\frac{g}{2}]}h(g-h)(12- k_g - \frac{2g+1}{h(g-h)})s_h \leq
(12- k_g - \frac{2g+1}{g^2 /4})\sum_{h=1}^{[\frac{g}{2}]}h(g-h)s_h .$$
By Corollary~\ref{hbound},
$$(12- k_g - \frac{2g+1}{g^2 /4})\sum_{h=1}^{[\frac{g}{2}]}h(g-h)s_h \leq
\frac{1}{4}[(3g+2)n-4(2g+1)](12-k_g -\frac{2g+1}{g^2 /4}).$$
Hence it suffices to find $k_g$ and $l_g$ such that 
$$ \frac{1}{4}[(3g+2)n-4(2g+1)](12-k_g -\frac{2g+1}{g^2 /4}) \leq 
[\frac{k_g}{4}g -(g-1)]n+ m_g.$$
First consider the inequality
$$ \frac{1}{4}(3g+2)(12-k_g -\frac{2g+1}{g^2 /4}) \leq
\frac{k_g}{4}g -(g-1).$$
Solving for $k_g$ we get,
$$k_g \geq 10 - \frac{6g+4}{g^2}.$$ 
Hence it suffices to find $l_g$ such that
$$-(2g+1)(12-k_g -\frac{2g+1}{g^2 /4}) 
\leq  m_g = -(g-1)(2g+1)k_g +(2g+1)l_g +8(g-1)(2g+1).$$
Solving for $l_g$ we get,
$$l_g \geq (k_g -8)(g-1) -12 + k_g +\frac{8g+4}{g^2} .$$
Therefore, if $g$ is even, then
$$k_g =10 - \frac{6g+4}{g^2} $$ 
and
$$l_g = 2g-10 + \frac{4g+4}{g^2}$$
will satify the requirements of the statement in the theorem.

Similarly, if $g$ is odd and $g \geq 3$, then one can calculate that
$$k_g = 10 - \frac{6g+4}{g^2 -1} $$ 
and
$$l_g =  2g-10 + \frac{8}{2g+1} + \frac{6}{(g^2 -1)(2g+1)}$$
will suffice. 

Moreover, we can choose $k_1 = l_1 = 0$ since $c_1^2=0$ 
for all genus one Lefschetz fibrations over $S^2$.

It is trivial to verify that $k_g \leq 10- \frac{6g+4}{g^2}$ and 
$l_g \leq 2g-10+ \frac{4g+4}{g^2}$ for all $g \geq 1$.

\end{proof}

{\Rems 
1. In {\cite{s2}}, Stipsicz proves that $0 \leq c_1^2 \leq 10 \chi$ for
the (relatively minimal) genus $g (\geq 1)$ Lefschetz fibrations over
closed surfaces of nonzero genus. However, his result does not extend to 
cover the fibrations over $S^2$. He also points out that 
$c_1^2 \leq 10 \chi +2g-1$ holds for any 
Lefschetz fibration.

2. In particular, Corollary~\ref{gengeo} yields $k_2=6$, $ l_2=-3 $,
$k_3=7.25$, $l_3=-2.75$, $k_4=8.25$ and $l_4=-0.75$. 

3. Combining our results for genus two and genus three Lefschetz 
fibrations, we have shown the following: The signature 
of a smooth 4-manifold which admits a hyperelliptic 
Lefschetz fibration of genus $g \leq 3$ over $S^2$ is negative. 

Hence we can formulate the following natural question: What 
is the minimal $g$ such that a smooth 4-manifold with 
positive signature admits a (hyperelliptic) 
Lefschetz fibration of genus $g$ over $S^2$ ? 

4. If one can improve the inequality $$\sigma(M) \leq n-s-4 $$ 
in Corollary~\ref{hbound} to $$\sigma(M) \leq n-s-4(g-1) $$
then one can prove that $\; c_1^2 \leq 10 \chi \;$ for all hyperelliptic 
Lefschetz fibrations over $S^2$.
}
\section{Examples}

\subsection{GENUS 1}

\indent

To illustrate how one can develop an algorithm using  
our main theorems to calculate the signatures of
smooth Lefschetz fibrations, we will give the details of
our computation 
to obtain the well-known result 
$\sigma(E(1))=-8$ for the elliptic surface $E(1)$.

\begin{figure}[h]
 \begin{center}
    \includegraphics{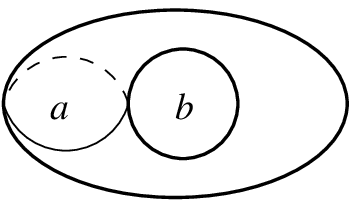}
  \caption{} \label{genus1}
   \end{center}
 \end{figure}

The global monodromy of $E(1)$ is given by the 
sequence $(\a , \b)^6$ of 12 Dehn
twists where $\a =D(a)$ and
$\b =D(b)$ denote the positive
Dehn twists about the curves $a$ and $b$, 
respectively (Figure 3). 

To build up $E(1)$, we start with a copy of 
$T^2 \times D^2$ and glue 2-handles along the vanishing cycles $a$ 
and $b$ in an alternating fashion. (We will use the letters 
$a$ and $b$ also to denote the homology classes of the curves $a$ and $b$, 
respectively.) 

Let $\phi$ denote the monodromy of the boundary of 
the 4-manifold before
we attach a 2-handle. 

Let $C =Ker(i_*)$ in the $\{[m^\prime],[l^\prime]\}$ plane, 
when we attach a 2-handle  (cf. 
Theorem~\ref{local}). Below we use Proposition~\ref{homology}
to compute $i_*[l^\prime]$ and $i_*[m^\prime]$.

$\phi =identity$, attach the first handle along $a$,

$i_*[l^\prime]=a$ and $i_*[m^\prime]=b-b^\prime=0$. 

$C =<[m^\prime]>$, $\sigma (id,a) =0$

$\phi = \a $, attach the second handle along $b$,

$i_*[l^\prime]=b$ and $i_*[m^\prime]=a^\prime -a=0$. 

$C =<[m^\prime]>$, $\sigma ( \a ,b) =0$

$\phi = \b \a$, attach the third handle along $a$, 

$i_*[l^\prime]=a$ and $i_*[m^\prime]=b-b^\prime=-a$. 

$C =<[m^\prime]+[l^\prime]>$, $\sigma (\b \a ,a)  =-1$

$\phi = \a \b \a$, attach the fourth handle along $b$, 

$i_*[l^\prime]=b$ and $i_*[m^\prime]=a^\prime -a=-2b$. 

$C =<[m^\prime]+2[l^\prime]>$, $\sigma (\a \b \a ,b) =-1$

$\phi = \b \a \b \a$, attach the fifth handle along $a$,

$i_*[l^\prime]=a$ and $i_*[m^\prime]=b-b^\prime=3b$. 

$C =<[m^\prime]+3[l^\prime]>$, $\sigma (\b \a \b \a ,a) =-1$

$\phi = \a \b \a \b \a$, attach the sixth handle along $b$,

$i_*[l^\prime]=b=0$ and $i_*[m^\prime]=a^\prime -a=-2a$. 

$C =<[l^\prime]>$, $\sigma (\a \b \a \b \a ,b) =-1$

$\phi = \b \a \b \a \b \a$, attach the seventh handle along $a$,

$i_*[l^\prime]=a=0$ and $i_*[m^\prime]=b-b^\prime=2b$. 

$C =<[l^\prime]>$, $\sigma ( \b \a \b \a \b \a ,a) =-1$

$\phi = \a \b \a \b \a \b \a$, attach the eighth handle along $b$,

$i_*[l^\prime]=b$ and $i_*[m^\prime]=a^\prime -a=4b$. 

$C =<-[m^\prime]+4[l^\prime]>$, $\sigma (\a \b \a \b \a \b \a ,b) =-1$
   
$\phi = \b \a \b \a \b \a \b \a$, attach the ninth handle along $a$,

$i_*[l^\prime]=a$ and $i_*[m^\prime]=b-b^\prime=3a$. 

$C =<-[m^\prime]+3[l^\prime]>$, $\sigma ( \b \a \b \a \b \a \b \a ,a)=-1$

$\phi = \a \b \a \b \a \b \a \b \a $, attach the tenth handle along $b$,

$i_*[l^\prime]=b$ and $i_*[m^\prime]=a^\prime -a=2b$. 

$C =<-[m^\prime]+2[l^\prime]>$, $\sigma ( \a \b \a \b \a \b \a \b \a ,b) =-1$

$\phi = \b \a \b \a \b \a \b \a \b \a $, attach the eleventh handle along $a$,

$i_*[l^\prime]=a$ and $i_*[m^\prime]=b-b^\prime=a$. 

$C =<-[m^\prime]+[l^\prime]>$, 
$\sigma ( \b \a \b \a \b \a \b \a \b \a ,a) =0$

$\phi = \a \b \a \b \a \b \a \b \a \b \a $, 
attach the twelfth handle along $b$,

$i_*[l^\prime]=b$ and $i_*[m^\prime]=a^\prime -a=b$. 

$C =<-[m^\prime]+[l^\prime]>$, 
$\sigma (\a \b \a \b \a \b \a \b \a \b \a ,b) =0$

Therefore by Theorem~\ref{burak1}

$ \sigma(E(1))= \sigma (id,a) + \sigma (\a ,b) + \sigma (\b \a ,a) +
\cdots + \sigma (\a \b \a \b \a \b \a \b \a \b \a ,b) =-8$. 

\subsection{GENUS 2}

\indent

We developed a Mathematica program to compute the signature 
of a 4-manifold which admits
a Lefschetz 
fibration over $D^2$ or $S^2$ whose global
monodromy is given by any finite
sequence of positive Dehn twists  
$D(c_1), D(c_2),..., D(c_5)$, where
$c_1,..., c_5$ are the curves indicated in Figure $4$.

\begin{figure}[h]
 \begin{center}
    \includegraphics{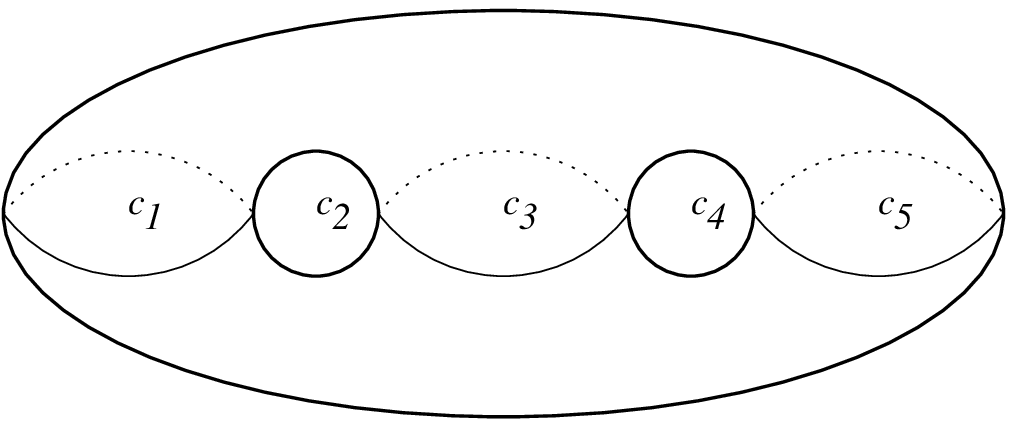}
  \caption{} \label{genus2}
   \end{center}
 \end{figure}

We computed the signatures of the following genus two Lefschetz 
fibrations 
of some closed $4$-manifolds,
using our Mathematica program (which is available upon request).

Let ${\z}_i$ denote $D(c_i)$, $1 \leq i \leq 5$.

$\sigma (({\z}_1,{\z}_2,{\z}_3,{\z}_4,{\z}_5,{\z}_5,{\z}_4,
{\z}_3,{\z}_2,{\z}_1)^2)=-12$,
 
$\sigma (({\z}_1,{\z}_2,{\z}_3,{\z}_4,{\z}_5)^6)=-18$,

$\sigma(({\z}_5,{\z}_1,{\z}_4,{\z}_2,{\z}_3,
{\z}_4,{\z}_2,{\z}_5,{\z}_1,{\z}_4,{\z}_2,{\z}_3,
{\z}_4,{\z}_2,{\z}_3)^2)=-18$,

$\sigma (({\z}_1,{\z}_2,{\z}_3,{\z}_4)^{10})=-24$.

{\Rems
1. One can indeed check these numbers using 
Matsumoto's local signature
formula. 

2. All of the examples above are fibrations with even number
of singular fibers. We want to point out that there are also examples
with odd number of singular fibers in a Lefschetz fibration.
Just take Matsumoto's example with $8$ singular fibers and 
replace one of the separating twists with the sequence of 
$12$ nonseparating twists $({\z}_1,{\z}_2)^6$ so that the resulting 
fibration will have $19$ singular fibers. If we replace the 
remaining separating twist with another sequence
of $12$ nonseparating twists then we get a fibration 
with $30$ singular fibers. It is a question whether 
the resulting manifold is 
diffeomorphic to $({\z}_1,{\z}_2,{\z}_3,{\z}_4,{\z}_5)^6$ or not.

We can generalize this to any genus two Lefschetz 
fibration over $S^2$ :

{\Prop \label{replace}  
Let $M$ be a $4$-manifold which admits a genus two
Lefschetz fibration over $S^2$.
Let $M^\prime$ be the
resulting manifold obtained by replacing a separating 
vanishing cycle in the global monodromy of this fibration 
by the sequence $({\z}_1,{\z}_2)^6$. Then
$$\sigma(M^\prime)=\sigma(M)-7$$ and 
$$e(M^\prime)=e(M)+11.$$}
\begin{proof}
We will give two different proofs, one using Matsumoto's local
signature formula, the other using the technique 
we developed in this paper.
Let $n$ and $s$ be the numbers of nonseparating and separating
vanishing cycles for $M$, respectively. 
Replacing one of the separating cycles by a product of $12$
nonseparating cycles we get $(n+12)$ nonseparating and 
$(s-1)$ separating cycles for $M^\prime$. Thus
$$\sigma(M^\prime)=(-3/5)(n+12)+(-1/5)(s-1)
=(-3/5)n+(-1/5)s-7=\sigma(M)-7.$$
Also using part $2$ in Theorem~\ref{local} when we remove a 
separating twist the signature will increase by $1$.
Gluing $({\z}_1,{\z}_2)^6$ we add $-8$ to the signature 
because the computation in genus two will be the same 
as in genus one ($\sigma((\a , \b)^6)=-8 $) using our 
technique. Hence we get the same result.     

Moreover $e(M^\prime)=e(M)+11 $  since we replace
a singular fiber with $12$ singular fibers.

\end{proof}

We also computed the signatures of the following genus two Lefschetz 
fibrations 
over $D^2$:

$\sigma(({\z}_1,{\z}_3,{\z}_5,{\z}_2,{\z}_4,
{\z}_2,{\z}_1,{\z}_1,{\z}_1,{\z}_3,{\z}_5,{\z}_3,
{\z}_5,{\z}_2,
{\z}_4,{\z}_5))=-10  $.

$\sigma(({\z}_1,{\z}_3,{\z}_5,{\z}_2,{\z}_4,
{\z}_2,{\z}_1,{\z}_1,{\z}_1,{\z}_3,{\z}_5,{\z}_3,
{\z}_5,{\z}_2,
{\z}_4,{\z}_5)^2)=-20  $.

$\sigma(({\z}_1,{\z}_3,{\z}_5,{\z}_2,{\z}_4,
{\z}_2,{\z}_1,{\z}_1,{\z}_1,{\z}_3,{\z}_5,{\z}_3,
{\z}_5,{\z}_2,
{\z}_4,{\z}_5)^3)=-30  $.

{\Rem We can not use the local signature
formula in order to check these results 
because it only works when the base space is $S^2$ (or any 
closed surface).}   

\subsection{GENUS 3}

\indent

We computed the signatures of the following genus 
three Lefschetz fibrations over $S^2$,
using our Mathematica program. 
(The program is available upon request.)

$\sigma (({\e}_1,{\e}_2,{\e}_3,{\e}_4,{\e}_5,{\e}_6)^{14})=-48$

$\sigma(({\e}_8,{\e}_9,{\e}_4,{\e}_3,{\e}_2,{\e}_1,{\e}_5,{\e}_4, 
{\e}_3,{\e}_2,{\e}_6,{\e}_5,{\e}_4,{\e}_3,({\e}_1,{\e}_2,{\e}_3,{\e}_4,
{\e}_5,{\e}_6)^{10}))=-42$.      

(This is a word in ${\CM}_3$ given in {\cite{f}}).

Here ${\e}_1, {\e}_2,..., {\e}_9 $ 
denote the positive Dehn twists about the curves
$d_1, d_2,..., d_9 $ indicated as in Figure $5$.

 \begin{figure}[h]
 \begin{center}
    \includegraphics{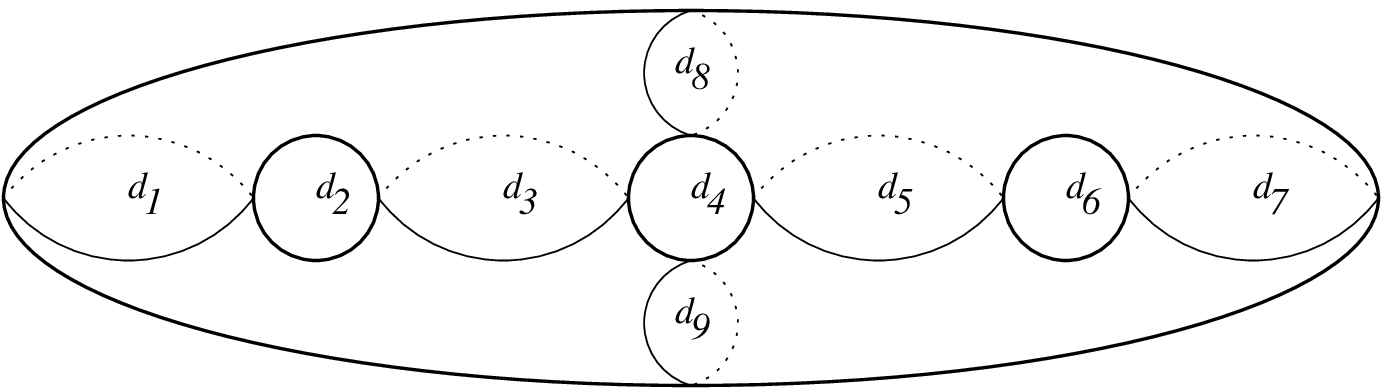}
  \caption{} \label{genus3}
  \end{center}
 \end{figure}

{\Rem The latter of these fibrations is not 
hyperelliptic since otherwise 
the local signature formula {\cite{e}} would yield 
$\sigma = 74 (-4/7)$ which is not an integer!}

We also computed the signatures of the following genus 
three Lefschetz fibrations over $D^2$:

$\sigma (({\e}_1,{\e}_3,{\e}_5,{\e}_7,{\e}_2,{\e}_4,{\e}_6))=-1 $

$\sigma (({\e}_1,{\e}_3,{\e}_5,{\e}_7,{\e}_2,{\e}_4,{\e}_6)^2)=-8 $

$\sigma (({\e}_1,{\e}_3,{\e}_5,{\e}_7,{\e}_2,{\e}_4,{\e}_6)^3)=-11 $

$\sigma (({\e}_1,{\e}_3,{\e}_5,{\e}_7,{\e}_2,{\e}_4,{\e}_6)^4)=-16 $

\section{Final Remark}

Given a product of positive Dehn twists in the mapping class group 
of a genus $ g$ surface, we can construct a symplectic 4-manifold which
admits a Lefschetz fibration over $D^2$, as we have studied in this paper.
A natural generalization is to allow negative Dehn twists also.
These fibrations are called {\em achiral} Lefschetz fibrations.
Our technique clearly extends to compute the signatures 
of these fibrations.

\end{document}